\documentclass{amsart}
\usepackage{amscd,amsthm}           

\newtheorem{theorem}{Theorem}[section]
\newtheorem{lemma}[theorem]{Lemma}
\newtheorem{proposition}[theorem]{Proposition}
\newtheorem{corollary}[theorem]{Corollary}
\newtheorem{conjecture}[theorem]{Conjecture}
\theoremstyle{definition}
\newtheorem{definition}[theorem]{Definition}

\DeclareMathOperator{\Spec}{Spec}
\DeclareMathOperator{\supp}{supp}
\DeclareMathOperator{\Hom}{Hom}
\DeclareMathOperator{\Ext}{Ext}

\DeclareMathOperator{\cok}{cok}

\newcommand{\cat}[1]{\mathcal{#1}}
\newcommand{\Ch}[1]{\text{Ch}(#1)}

\newcommand{\ideal}[1]{\mathfrak{#1}}

\newcommand{\Z}{\mathbb{Z}}

\newcommand{\Mod}{\text{-mod}}

\newcommand{\mathcolon}{\colon\,}

\newcommand{\wide}{\cat{L}_{\textup{wide}}}
\newcommand{\serre}{\cat{L}_{\textup{Serre}}}
\newcommand{\tors}{\cat{L}_{\textup{tors}}}
\newcommand{\thick}{\cat{L}_{\textup{thick}}}
\newcommand{\localizing}{\cat{L}_{\textup{wide}}^{\oplus }}
\newcommand{\mloc}{\thick ^{\oplus }}
\begin{document}
 
\title{Classifying subcategories of modules} 

\date{\today}

\author{Mark Hovey}
\address{Department of Mathematics \\ Wesleyan University
\\ Middletown, CT 06459}
\email{hovey@member.ams.org}




\maketitle

\section*{Introduction}

A basic problem in mathematics is to classify all objects one is
studying up to isomorphism.  A lesson this author learned from stable
homotopy theory~\cite{hopkins-smith} is that while this is almost always
impossible, it is sometimes possible, and very useful, to classify
collections of objects, or certain full subcategories of the category
one is working in.  In particular, if $R$ is a commutative ring, the
thick subcategories of small objects in the derived category
$\cat{D}(R)$ have been classified.  Recall that a thick subcategory is a
triangulated subcategory closed under summands.  Thick subcategories
correspond to unions of subsets of $\Spec R$ of the form $V(\ideal{a})$,
where $\ideal{a}$ is a finitely generated ideal of $R$.  In particular,
when $R$ is Noetherian, they correspond to arbitrary unions of closed
sets of $\Spec R$.  This line of research was initiated by
Hopkins~\cite{hopkins-global-methods}, where he wrote a down a false
proof of this classification.  Neeman~\cite{neeman-derived} later
corrected Hopkins' proof in the Noetherian case, and
Thomason~\cite{thomason-tri} generalized the result to arbitrary
commutative rings (and, in fact, to quasi-compact, quasi-separated
schemes).

The author has long thought that the analogous classification in the
ostensibly simpler category of $R$-modules is the classification of
torsion theories when $R$ is a Noetherian and commutative
ring~\cite[Section~VI.6]{stenstrom}.  After all, these too correspond to
arbitrary unions of closed sets in $\Spec R$.  However, we show in this
paper that the analog of a thick subcategory in $\cat{D}(R)$ is not a
torsion theory of $R$-modules, but just an Abelian subcategory of
$R$-modules closed under extensions.  We call this a wide subcategory.
We use the classification of thick subcategories mentioned above to give
a classification of wide subcategories of finitely presented modules
over a large class of commutative coherent rings.  To be precise, our
classification works for quotients of regular commutative coherent rings
by finitely generated ideals.  Recall that a coherent ring is regular if
every finitely generated ideal has finite projective dimension.  Thus
our classification includes, for example, the polynomial ring on
countably many variables over a principal ideal domain, and all finitely
generated algebras over a field.  A corollary of our result is that
if $R$ is Noetherian as well, then every wide subcategory of finitely
presented $R$-modules is in fact the collection of all finitely
presented modules in a torsion theory.  It is interesting that we have
no direct proof of this fact, but must resort to the rather difficult
classification of thick subcategories in the derived category.  

One can also attempt to classify thick subcategories closed under
arbitrary coproducts, or arbitrary products.  These are called
localizing and colocalizing subcategories, respectively.  For
$\cat{D}(R)$ when $R$ is a Noetherian commutative ring, they were
classified by Neeman~\cite{neeman-derived}, and correspond to arbitrary
subsets of $\Spec R$.  We give an analogous classification of wide
subcategories closed under arbitrary coproducts in the Noetherian case.
Once again, the proof of this relies on comparison with the derived
category.  

\section{Wide subcategories}

Suppose $R$ is a ring.  In this section, we define wide subcategories of
$R$-modules and construct an adjunction to thick subcategories of
$\cat{D}(R)$.  The natural domain of this adjunction is actually $\wide
(R)$, the lattice of wide subcategories of $\cat{C}_{0}$, the wide
subcategory generated by $R$.  When $R$ is coherent, we identify
$\cat{C}_{0}$ with the finitely presented modules, but we do not know
what it is in general.

We recall that a thick subcategory of a triangulated category like
$\cat{D}(R)$ is a full triangulated subcategory closed under retracts
(summands).  This means, in particular, that if we have an exact
triangle $X\xrightarrow{}Y\xrightarrow{}Z\xrightarrow{}\Sigma X$ and two
out of three of $X$, $Y$, and $Z$ are in the thick subcategory, so is
the third.

The analogous definition for subcategories of an Abelian category is the
following.  

\begin{definition}\label{defn-zaftig}
A full subcategory $\cat{C}$ of $R\Mod $, or any Abelian category, is
called \emph{wide} if it is Abelian and closed under extensions.
\end{definition}

When we say a full subcategory $\cat{C}$ is Abelian, we mean that if
$f\mathcolon M\xrightarrow{}N$ is a map of $\cat{C}$, then the kernel
and cokernel of $f$ are in $\cat{C}$.  Thus a wide subcategory
$\cat{C}$ need not be closed under arbitrary subobjects or quotient
objects.  However, $\cat{C}$ is automatically closed under summands.
Indeed, if $M\cong N\oplus P$ and $M\in \cat{C}$, then $N$ is the kernel
of the self-map of $M$ that takes $(n,p)$ to $(0,p)$.  Thus $N\in
\cat{C}$.  In particular, $\cat{C}$ is replete, in the sense that
anything isomorphic to something in $\cat{C}$ is itself in $\cat{C}$.

Torsion theories and Serre classes are closely related to wide
subcategories.  Recall that a Serre class is just a wide subcategory
closed under arbitrary subobjects, and hence arbitrary quotient objects.
Similarly, a (hereditary) torsion theory is a Serre class closed under
arbitrary direct sums.  In particular, the empty subcategory, the $0$
subcategory, and the entire category of $R$-modules are torsion
theories, and so wide subcategories.  The category of all $R$-modules
of cardinality $\leq \kappa $ for some infinite cardinal $\kappa $ is a
Serre class (but not a torsion theory), and hence a wide
subcategory.  The category of finite-dimensional rational vector spaces,
as a subcategory of the category of abelian groups, is an example of a
wide subcategory that is not a Serre class.  The thick subcategories
studied in~\cite{hovey-palmieri-quillen, hovey-palmieri-galois} are, on
the other hand, more general than wide subcategories.  

Note that the collection of all wide subcategories of $R\Mod $ forms
a complete lattice (though it is not a set).  Indeed, the join of a
collection of wide subcategories is the wide subcategory generated
by them all, and the meet of a collection of wide subcategories is
their intersection.  

The following proposition shows that wide subcategories are the
analogue of thick subcategories.  

\begin{proposition}\label{prop-zaftig-gives-thick}
Suppose $\cat{C}$ is a wide subcategory of $R\Mod $.  Define
$f(\cat{C})$ to be the collection of all small objects $X\in \cat{D}(R)$
such that $H_{n}X\in \cat{C}$ for all $n\in \Z $.  Then $f(\cat{C})$ is
a thick subcategory.  
\end{proposition}

Note that the collection of all thick subcategories is also a complete
lattice, and the map $f$ is clearly order-preserving.  

\begin{proof}
Since wide subcategories are closed under summands, $f(\cat{C})$ is
closed under retracts.  It is clear that $X\in f(\cat{C})$ if and only
if $\Sigma X\in f(\cat{C})$.  It remains to show that, if we have an
exact triangle $X\xrightarrow{}Y\xrightarrow{}Z\xrightarrow{}\Sigma X$
and $X,Z\in f(\cat{C})$, then $Y\in f(\cat{C})$.  We have a short exact
sequence
\[
0\xrightarrow{} A\xrightarrow{}H_{n}Y \xrightarrow{} B \xrightarrow{}0
\]
where $A$ is the cokernel of the map $H_{n+1}Z \xrightarrow{}H_{n}X$,
and $B$ is the kernel of the map $H_{n}Z\xrightarrow{}H_{n-1}X$.  Hence
$A$ and $B$ are in $\cat{C}$, and so $H_{n}Y\in \cat{C}$ as well.  Thus
$Y\in f(\cat{C})$.  
\end{proof}

Note that Proposition~\ref{prop-zaftig-gives-thick} remains true if
$R\Mod $ is replaced by any Abelian category, or indeed, if $\cat{D}(R)$
is replaced by a stable homotopy category~\cite{hovey-axiomatic} and $R$
is replaced by the homotopy of the sphere in that category.  

Proposition~\ref{prop-zaftig-gives-thick} implies that the homology of a
small object in $\cat{D}(R)$ must lie in the wide subcategory
generated by $R$.  

\begin{corollary}\label{cor-small}
Let $\cat{C}_{0}$ be the wide subcategory generated by $R$.  If $X$
is a small object of $\cat{D}(R)$, then $H_{n}X\in \cat{C}_{0}$ for all
$n$ and $H_{n}X=0$ for all but finitely many $n$.  
\end{corollary}

\begin{proof}
By hypothesis, the complex $S^{0}$ consisting of $R$ concentrated in
degree $0$ is in $f(\cat{C}_{0})$.  Therefore the thick subcategory
$\cat{D}$ generated by $S^{0}$ is contained in $f(\cat{C}_{0})$.  But
$\cat{D}$ is precisely the small objects in $\cat{D}(R)$.  (This is
proved in~\cite[Corollary~2.3.12]{hovey-axiomatic} for commutative $R$,
but the proof does not require commutativity).  Hence $H_{n}X\in
\cat{C}_{0}$ for all small objects $X$ and all $n$.  It remains to
prove that $H_{n}X=0$ for all but finitely many $n$ if $X$ is small.
This is proved analogously; the collection of all such $X$ is a thick
subcategory containing $S^{0}$.
\end{proof}

This corollary tells us that the proper domain of $f$ is $\wide (R)$,
the lattice of wide subcategories of $\cat{C}_{0}$.  We would like $f$
to define an isomorphism $f\mathcolon \wide (R)\xrightarrow{}\thick
(\cat{D}(R))$, where $\thick (\cat{D}(R))$ is the lattice of thick
subcategories of small objects in $\cat{D}(R)$.  We now construct the
only possible inverse to $f$.

Given a thick subcategory $\cat{D}\in \thick (\cat{D}(R))$, we define
$g(\cat{D})$ to be the wide subcategory generated by $\{H_{n}X \}$,
where $X$ runs though objects of $\cat{D}$ and $n$ runs through $\Z $.
By Corollary~\ref{cor-small}, $g(\cat{D})\in \wide (R)$.  Also,
$g$ is obviously order-preserving.  We also point out that, like $f$,
$g$ can be defined in considerably greater generality.

We then have the following proposition.  

\begin{proposition}\label{prop-adjoint}
The lattice homomorphism $g$ is left adjoint to $f$.  That is, for
$\cat{C}\in \wide (R)$ and $\cat{D}\in \thick (\cat{D}(R))$, we have
$g(\cat{D})\subseteq \cat{C}$ if and only if $\cat{D}\subseteq
f(\cat{C})$.
\end{proposition}

\begin{proof}
Suppose first that $g(\cat{D})\subseteq \cat{C}$.  This means that for
every $X\in \cat{D}$, we have $H_{n}X\in \cat{C}$ for all $n$.  Hence
$X\in f(\cat{C})$.  Thus $\cat{D}\subseteq f(\cat{C})$.  Conversely, if
$\cat{D}\subseteq f(\cat{C})$, then for every $X\in \cat{D}$ we have
$H_{n}X\in \cat{C}$ for all $n$.  Thus $g(\cat{D})\subseteq \cat{C}$. 
\end{proof}

\begin{corollary}\label{cor-adjoint}
Suppose $R$ is a ring.  If $\cat{C}\in \wide (R)$, then $gf(\cat{C})$ is
the smallest wide subcategory $\cat{C}'$ such that
$f(\cat{C}')=f(\cat{C})$.  Similarly, if $\cat{D}\in \thick
(\cat{D}(R))$, then $fg(\cat{D})$ is the largest thick subcategory
$\cat{D}'$ such that $g(\cat{D}')=g(\cat{D})$.
\end{corollary}

\begin{proof}
This corollary is true for any adjunction between partially ordered
sets.  For example, if $f(\cat{C}')= f(\cat{C})$, then
$gf(\cat{C}')=gf(\cat{C})$.  But $gf(\cat{C}')\subseteq \cat{C}'$, so
$gf(\cat{C})\subseteq \cat{C}'$.  Furthermore, combining the counit and
unit of the adjunction shows that $fgf(\cat{C})$ is contained in and
contains $f(\cat{C})$.  The other half is similar.  
\end{proof}

It follows from this corollary that $f$ is injective if and only
$gf(\cat{C})=\cat{C}$ for all $\cat{C}\in \wide (R)$ and that $f$ is
surjective if and only if $fg(\cat{D})=\cat{D}$ for all $\cat{D}\in
\thick (\cat{D}(R))$.  

In order to investigate these questions, it would be a great help to
understand $\cat{C}_{0}$, the wide subcategory generated by $R$.  We
know very little about this in general, except that $\cat{C}_{0}$
obviously contains all finitely presented modules and all finitely
generated ideals of $R$.  We also point out that $\cat{C}_{0}$ is
contained in the wide subcategory consisting of all modules of
cardinality $\leq \kappa $, where $\kappa $ is the larger of $\omega $
and the cardinality of $R$.  In particular, $\cat{C}_{0}$ has a small
skeleton, and so there is only a set of wide subcategories of
$\cat{C}_{0}$.

The only case where we can identify $\cat{C}_{0}$ is when $R$ is a
coherent ring.  A brief description of coherent rings can be found
in~\cite[Section~I.13]{stenstrom}; an excellent reference for deeper
study is~\cite{glaz}.  

\begin{lemma}\label{lem-coherent}
A ring $R$ is coherent if and only if the wide subcategory
$\cat{C}_{0}$ generated by $R$ consists of the finitely presented
modules.
\end{lemma}

\begin{proof}
Suppose first that $\cat{C}_{0}$ is the collection of finitely presented
modules.  Suppose $\ideal{a}$ is a finitely generated left ideal of
$R$.  Then $R/\ideal{a}$ is a finitely presented module, so $\ideal{a}$,
as the kernel of the map $R\xrightarrow{}R/\ideal{a}$, is in
$\cat{C}_{0}$.  Hence $\ideal{a}$ is finitely presented, and so $R$ is
coherent.  

The collection of finitely presented modules over any ring is clearly
closed under cokernels and is also closed under extensions.  Indeed,
suppose we have a short exact sequence
\[
0\xrightarrow{}M'\xrightarrow{}M\xrightarrow{}M''\xrightarrow{}0
\]
where $M'$ and $M''$ are finitely presented (in fact, we need only
assume $M'$ is finitely generated).  Choose a finitely generated
projective $P$ and a surjection $P\xrightarrow{}M''$.  We can lift this
to a map $P\xrightarrow{}M$.  Then we get a surjection $M'\bigoplus
P\xrightarrow{}M$, as is well-known.  Furthermore, the kernel of this
surjection is the same as the kernel of $P\xrightarrow{}M''$, which is
finitely generated since $M''$ is finitely presented.  Hence $M$ is
finitely presented.  

Now suppose that $R$ is coherent.  We show that the kernel of a map
$f\mathcolon M\xrightarrow{}N$ of finitely presented modules is finitely
presented.  The point is that the image of $f$ is a finitely generated
submodule of the finitely presented module $N$.  Because the ring is
coherent, this means that the image of $f$ is finitely presented.  The
kernel of $f$ is therefore finitely generated, but it is a submodule of
the finitely presented module $M$, so it is finitely presented, using
coherence again.
\end{proof}

Noetherian rings can be characterized in a similar manner as rings in
which $\cat{C}_{0}$ is the collection of finitely generated modules.  

\section{Surjectivity of the adjunction}\label{sec-surj}

The goal of this section is to show that the map $f\mathcolon \wide
(R)\xrightarrow{}\thick (\cat{D}(R))$ is surjective for all commutative
rings $R$.  This is a corollary of Thomason's classification of thick
subcategories in $\cat{D}(R)$.

Suppose $R$ is a commutative ring.  Denote by $J(\Spec R)\subseteq
2^{\Spec R}$ the collection of order ideals in $\Spec R$, so that $S\in
J(\Spec R)$ if and only if $\ideal{p}\in S$ and $\ideal{q}\subseteq
\ideal{p}$ implies that $\ideal{q}\in S$.  Note that an open set in the
Zariski topology of $\Spec R$ is in $J(\Spec R)$, so an arbitrary
intersection of open sets is in $J(\Spec R)$.  Also, note that $J(\Spec
R)$ is a complete distributive lattice.  

We will construct a chain of maps 
\[
J(\Spec R)^{\textup{op}} \xrightarrow{i} \tors (R) \xrightarrow{j} \serre (R)
\xrightarrow{\alpha } \wide (R) \xrightarrow{f} \thick (\cat{D}(R)),
\]
each of which is a right adjoint.  We have of course already constructed
$f$.  

To construct $i$, note that $\tors (R)$ denotes the lattice of all
torsion theories of $R$-modules.  Recall that a torsion theory is a wide
subcategory closed under arbitrary submodules and arbitrary direct sums.
The map $i$ is defined by 
\[
i(S)=\{M|M_{(\ideal{p})}=0 \text{ for all }
\ideal{p} \in S\}.
\]
Its right adjoint $r$ has
\[
r(\cat{T})= \{\ideal{p}| M_{(\ideal{p})}=0 \text{ for all } M\in \cat{T}
\} = \bigcap _{M\in \cat{T}} (\Spec R \setminus \supp M). 
\]
One can check that $ri(S)=S$, so that $i$ is an
embedding.  In case $R$ is a Noetherian commutative ring, $i$ is an
isomorphism~\cite[Section~VI.6]{stenstrom}.

To construct $j$, let $\cat{S}_{0}$ denote the Serre class generated by
$R$.  Recall that a Serre class is a wide subcategory closed under
arbitrary subobjects.  If $R$ is Noetherian, then $\cat{S}_{0}$ is the
finitely generated $R$-modules, but in general it will be larger than
this.  The symbol $\serre (R)$ denotes the lattice of Serre subclasses
of $\cat{S}_{0}$.  The map $j$ takes a torsion theory $\cat{T}$ to its
intersection with $\cat{S}_{0}$.  Its left adjoint $s$ takes a Serre
subclass of $\cat{S}_{0}$ to the torsion theory it generates.  Since a
torsion theory is determined by the finitely generated modules in it
(since it is closed under direct limits), the composite $sj$ is the
identity.  Thus $j$ is also an embedding, for an arbitrary ring $R$.
When $R$ is Noetherian, $j$ is an isomorphism.  Indeed, in this case,
the collection of modules all of whose finitely generated submodules lie
in a Serre subclass $\cat{S}$ of finitely generated modules is a torsion
theory, and is therefore $s(\cat{S})$.  Hence $js$ is the identity as
well.

The map $\alpha $ takes a Serre class to its intersection with
$\cat{C}_{0}$.  Its adjoint $\beta $ takes a wide subcategory to the
Serre class it generates.  When $R$ is Noetherian, $\cat{C}_{0}$ and
$\cat{S}_{0}$ coincide, so one can easily see that $\beta \alpha $ is
the identity, so that $\alpha $ is injective.  However, $\alpha $ will not
be injective in general, as we will see below.  

Note that the composite $f\alpha ji$ takes $S\in J(\Spec R)$ to the
collection of all small objects $X$ in $\cat{D}(R)$ such that
$(H_{n}X)_{(\ideal{p})}=0$ for all $\ideal{p}\in S$ and all $n$.  Since
$H_{n}(X_{(\ideal{p})})\cong (H_{n}X)_{(\ideal{p})}$, this is the same
as the collection of all small $X$ such that $X_{(\ideal{p})}=0$ for all
$\ideal{p}\in S$.

The following theorem is the main result of~\cite{thomason-tri}.  To
describe it, recall that the open subsets of the Zariski topology on
$\Spec R$, where $R$ is commutative, are the sets $D(\ideal{a})$ where
$\ideal{a}$ is an ideal of $R$ and $D(\ideal{a})$ consists of all primes
that do not contain $\ideal{a}$.  The open set $D(\ideal{a})$ is
quasi-compact if and only if $D(\ideal{a})=D(\ideal{b})$ for some
finitely generated ideal $\ideal{b}$ of $R$.  This fact is well-known in
algebraic geometry, and can be deduced from the argument on the top of
p. 72 in~\cite{hartshorne}.  Now we let $\widetilde{J}(\Spec R)$ denote
the sublattice of $J(\Spec R)$ consisting of arbitrary intersections of
quasi-compact open sets.  

\begin{theorem}[Thomason's theorem]\label{thm-thomason}
Let $R$ be a commutative ring.  Let $h$ denote the restriction of
$f\alpha ji$ to $\widetilde{J}(\Spec R)$.  Then $h\mathcolon
\widetilde{J}(\Spec R)^{\textup{op}}\xrightarrow{}\thick (\cat{D}(R))$
is an isomorphism.
\end{theorem}

The following corollary is immediate, since $f\alpha ji$ is surjective.  

\begin{corollary}\label{cor-thomason}
Let $R$ be a commutative ring.  Then the map $f\mathcolon \wide
(R)\xrightarrow{}\thick (\cat{D}(R))$ is surjective.  In particular, for
any thick subcategory $\cat{D}$, we have $fg(\cat{D})=\cat{D}$.
\end{corollary}

Note that, since $i$ and $j$ are injective for all rings $R$, torsion
theories and Serre classes cannot classify thick subcategories of
$\cat{D}(R)$ in general.  There are torsion theories and Serre classes
of $R$-modules that do not correspond to any thick subcategory of small
objects in $\cat{D}(R)$.  When $R$ is Noetherian, we have
$\widetilde{J}(\Spec R)=J(\Spec R)$, so torsion theories and Serre
classes do correspond to thick subcategories, but this will not be true
in general.  

\section{Regular coherent rings}\label{sec-injective}

The goal of this section is to show that the map $f\mathcolon \wide
(R)\xrightarrow{}\thick (\cat{D}(R))$ is an isomorphism when $R$ is a
regular coherent commutative ring.  Regularity means that every finitely
presented module has finite projective dimension;
see~\cite[Section~6.2]{glaz} for many results about regular coherent
rings.  An example of a regular coherent ring that is not Noetherian is
the polynomial ring on infinitely many variables over a principal ideal
domain.  

As we have already seen, $f$ is injective if and only
$gf(\cat{C})=\cat{C}$ for all wide subcategories of finitely
presented modules (when $R$ is coherent).  We start out by proving that
$gf(\cat{C}_{0})=\cat{C}_{0}$ when $R$ is coherent.

\begin{proposition}\label{prop-fin-pres}
Suppose $R$ is a ring, and $M$ is a finitely presented $R$-module.  Then
there is a small object $X\in \cat{D}(R)$ such that $H_{0}X\cong M$.  
\end{proposition}

\begin{proof}
Write $M$ as the cokernel of a map $f\mathcolon
R^{m}\xrightarrow{}R^{n}$.  Recall that, given a module $N$, $S^{0}N$
denotes the complex that is $N$ concentrated in degree $0$.  Define $X$
to be the cofiber of the induced map
$S^{0}R^{m}\xrightarrow{}S^{0}R^{n}$.  Then $X$ is small and $H_{0}X=M$
(and $H_{1}X$ is the kernel of $f$).
\end{proof}

\begin{corollary}\label{cor-fin-pres}
Suppose $R$ is a coherent ring, and $\cat{C}_{0}$ is the subcategory of
finitely presented modules.  Then $gf(\cat{C}_{0})=\cat{C}_{0}$.  
\end{corollary}

In order to prove that $gf(\cat{C})=\cat{C}$ in general, however, given
a finitely presented module $M$, we would have to find a complex $X$
that is small in $\cat{D}(R)$ such that $H_{0}X\cong M$ and each
$H_{n}X$ is in the wide subcategory generated by $M$.  The obvious
choice is $S^{0}M$, the complex consisting of $M$ concentrated in degree
$0$.  However, $S^{0}M$ cannot be small in $\cat{D}(R)$ unless $M$ has
finite projective dimension, as we show in the following lemma. 

\begin{lemma}\label{lem-fin-dim}
Suppose $R$ is a ring and $M$ is an $R$-module.  If the complex $S^{0}M$
is small in $\cat{D}(R)$, then $M$ has finite projective dimension.  
\end{lemma}

\begin{proof}
Define an object $X$ of $\cat{D}(R)$ to have finite projective dimension
if there is an $i$ such that $H_{j}F(X,S^{0}N)=0$ for all $R$-modules
$N$ and all $j$ with $|j|>i$.  Here $F(X,S^{0}N)$ is the function
complex $\Hom_{R} (QX,N)$ in $\cat{D}(\Z )$ obtained by replacing $X$ by
a cofibrant chain complex $QX$ quasi-isomorphic to it.  (In the
terminology of~\cite[Chapter~4]{hovey-model}, the model category
$\Ch{R}$ of chain complexes over $R$ with the projective model structure
is a $\Ch{\Z }$-model category, and we are using that structure).  If
$X=S^{0}M$, then $H_{i}F(S^{0}M,S^{0}N)=\Ext ^{-i}(M,N)$, so $S^{0}M$
has finite projective dimension if and only if $M$ does.  It is easy to
see that complexes with finite projective dimension form a thick
subcategory containing $R$.  Therefore every small object of
$\cat{D}(R)$ has finite projective dimension.
\end{proof}

Conversely, we have the following proposition.  

\begin{proposition}\label{prop-fin-proj-dim}
Suppose $R$ is a coherent ring and $M$ is a finitely presented module of
finite projective dimension.  Then $S^{0}M$ is small in $\cat{D}(R)$.  
\end{proposition}

\begin{proof}
It may be possible to give a direct proof of this, but we prefer to use
model categories.  Theorem~7.4.3 of~\cite{hovey-model} asserts that any
cofibrant complex $A$ that is small in the category $\Ch{R}$ of chain
complexes and chain maps, in the sense that $\Ch{R}(A,-)$ commutes with
direct limits, will be small in $\cat{D}(R)$.  Of course, $S^{0}M$ is
small in $\Ch{R}$, but it will not be cofibrant.  To make it cofibrant,
we need to replace $M$ by a projective resolution.  Since $M$ is
finitely presented and the ring $R$ is coherent, each term $P_{i}$ in a
projective resolution for $M$ will be finitely generated.  Since $M$ has
finite projective dimension, the resolution $P_{*}$ is finite.  Hence
$P_{*}$ is small in $\Ch{R}$, and so also in $\cat{D}(R)$.  Since
$P_{*}$ is isomorphic to $S^{0}M$ in $\cat{D}(R)$, the result follows.  
\end{proof}

This proposition leads immediately to the following theorem.  

\begin{theorem}\label{thm-thick-gives-zaftig}
Suppose $R$ is a regular coherent ring, and $\cat{C}$ is a wide
subcategory of finitely presented $R$-modules.  Then
$gf(\cat{C})=\cat{C}$.
\end{theorem}

\begin{proof}
We have already seen that $gf(\cat{C})\subseteq \cat{C}$.  Suppose $M\in
\cat{C}$.  Then $S^{0}M$ is small by
Proposition~\ref{prop-fin-proj-dim}, so clearly $S^{0}M\in f(\cat{C})$.
Thus $M\in gf(\cat{C})$.  
\end{proof}

The author believes that this theorem should hold without the regularity
hypothesis, though obviously a cleverer proof is required.
Theorem~\ref{thm-thick-gives-zaftig} and Corollary~\ref{cor-thomason}
lead immediately to the following classification theorem.

\begin{theorem}\label{thm-classification}
Suppose $R$ is a regular commutative coherent ring. Then the map
$f\mathcolon \wide (R)\xrightarrow{} \thick (\cat{D}(R))$ is an
isomorphism.  Hence the restriction of $\alpha ji$ defines an
isomorphism $\widetilde{h}\mathcolon \widetilde{J}(\Spec
R)\xrightarrow{}\wide (R)$ as well.
\end{theorem}

\begin{corollary}\label{cor-Serre}
Suppose $R$ is a regular commutative coherent ring, $\cat{C}$ is a wide
subcategory of finitely presented $R$-modules, and $M\in \cat{C}$.  If
$N$ is a finitely presented submodule or quotient module of $M$, then
$N\in \cat{C}$.  In particular, if $R$ is also Noetherian, every wide
subcategory of finitely generated modules is a Serre class.
\end{corollary}

Indeed, the first statement is obviously true for any wide subcategory
coming from $\widetilde{J}(\Spec R)$.  

\section{Quotients of regular coherent rings}\label{sec-ideal}

The goal of this section is to understand the relationship between wide
subcategories of $R$-modules and wide subcategories of
$R/\ideal{a}$-modules, where $\ideal{a}$ is a two-sided ideal of $R$.
This will allow us to extend Theorem~\ref{thm-classification} to
quotients of regular commutative coherent rings by finitely generated
ideals. 

Given $\cat{C}\in \wide (R/\ideal{a})$, we can think of $\cat{C}$ as a
full subcategory of $R$-modules where $\ideal{a}$ happens to act
trivially.  As such, it will be closed under kernels and cokernels, but
not extensions.  Define $u(\cat{C})$ to be the wide subcategory of
$R$-modules generated by $\cat{C}$.  In order to be sure that
$u(\cat{C})$ is contained in $\cat{C}_{0}(R)$, we need to make sure that
$R/\ideal{a}\in \cat{C}_{0}(R)$.  The easiest way to be certain of this
is if $\ideal{a}$ is finitely generated as a left ideal; under this
assumption, we have just defined a map $u\mathcolon \wide
(R/\ideal{a})\xrightarrow{}\wide (R)$.

As usual, this map has a left adjoint $v$.  Given $\cat{D}\in
\wide (R)$, we define $v(\cat{D})$ to be the collection
of all $M\in \cat{D}$ such that $\ideal{a}$ acts trivially on $M$, so
that $M$ is naturally an $R/\ideal{a}$-module.   Then $v(\cat{D})$ is
a wide subcategory.  It is not clear that $v(\cat{D})\subseteq
\cat{C}_{0}(R/\ideal{a})$ in general.  However, if $R$ is coherent, then
any $M$ in $g(\cat{D})$ will be finitely presented as an $R$-module, and
so finitely presented as an $R/\ideal{a}$-module.  

Altogether then, we have the following lemma, whose proof we leave to
the reader.  

\begin{lemma}\label{lem-mod-adjoint}
Suppose $R$ is a coherent ring and $\ideal{a}$ is a two-sided ideal of
$R$ that is finitely generated as a left ideal.  Then the map
$u\mathcolon \wide (R/\ideal{a})\xrightarrow{}\wide (R)$ constructed
above is right adjoint to the map $v$ constructed above.
\end{lemma}

We claim that $vu(\cat{C})=\cat{C}$ for all $\cat{C}\in
\wide (R/\ideal{a})$, so that $u$ is in fact an
embedding.  To see this, we need a description of $u(\cat{C})$, or, more
generally, a description of the wide subcategory generated by a full
subcategory $\cat{D}$ that is already closed under kernels and
cokernels.  It is clear that this wide subcategory will have to
contain all extensions of $\cat{D}$, so let $\cat{D}_{1}$ denote the
full subcategory consisting of all extensions of $\cat{D}$.  An
object $M$ is in $\cat{D}_{1}$ if and only if there is a short exact
sequence 
\[
0 \xrightarrow{} M' \xrightarrow{} M \xrightarrow{} M'' \xrightarrow{} 0
\]
where $M'$ and $M''$ are in $\cat{D}$.  Since $\cat{D}$ is closed under
kernels and cokernels, $0\in \cat{D}$ (unless $\cat{D}$ is empty), so
that $\cat{D}_{1}\supseteq \cat{D}$.  

We claim that $\cat{D}_{1}$ is still closed under kernels and
cokernels.  We will prove this after the following lemma. 

\begin{lemma}\label{lem-extensions}
Suppose $\cat{A}$ is an abelian category and the full subcategory
$\cat{D}$ of $\cat{A}$ is closed under kernels and cokernels.  Let
$\cat{D}_{1}$ be the full subcategory consisting of extensions of
$\cat{D}$.  Suppose $M\in \cat{D}$, $N\in \cat{D}_{1}$, and $f\mathcolon
M\xrightarrow{}N$ is a map.  Then $\ker f\in \cat{D}$ and $\cok f\in
\cat{D}_{1}$.
\end{lemma}

\begin{proof}
Because $N\in \cat{D}_{1}$, we can construct the commutative diagram
below,
\[
\begin{CD}
0 @>>> 0 @>>> M @= M @>>> 0 \\
@. @VVV @VfVV @Vp\circ fVV @. \\
0 @>>> N' @>>> N @>>p> N'' @>>> 0
\end{CD}
\]
where the rows are exact and $N',N''\in \cat{D}$.  The snake lemma then
gives us an exact sequence 
\[
0 \xrightarrow{} \ker f \xrightarrow{} \ker (p\circ f)
\xrightarrow{\partial } N' \xrightarrow{} \cok f \xrightarrow{} \cok
(p\circ f) \xrightarrow{} 0.
\]
Since $\ker (p\circ f)$ and $N'$ are in $\cat{D}$, we find that $\ker f
= \ker \partial $ is in $\cat{D}$.  Similary, we find that $\cok f$ is
an extension of $\cok \partial $ and $\cok (p\circ f)$, so $\cok f \in
\cat{D}_{1}$. 
\end{proof}

\begin{proposition}\label{prop-extensions}
Suppose $\cat{A}$ is an abelian category, and $\cat{D}$ is a full
subcategory of $\cat{A}$ closed under kernels and cokernels.  Let
$\cat{D}_{1}$ be the full subcategory consisting of extensions of
$\cat{D}$.  Then $\cat{D}_{1}$ is also closed under kernels and
cokernels.
\end{proposition}

\begin{proof}
Suppose $f\mathcolon M\xrightarrow{}N$ is a map of $\cat{D}_{1}$.  Then
we have the commutative diagram below, 
\[
\begin{CD}
0 @>>> M' @>i>> M @>>> M'' @>>> 0 \\
@. @Vf\circ iVV @VfVV @VVV @. \\
0 @>>> N @= N @>>> 0 @>>> 0
\end{CD}
\]
where the rows are exact and $M',M''\in \cat{D}$.  The snake lemma gives
an exact sequence 
\[
0 \xrightarrow{} \ker (f\circ i) \xrightarrow{} \ker f \xrightarrow{}
M'' \xrightarrow{\partial } \cok (f\circ i) \xrightarrow{} \cok f
\xrightarrow{} 0.
\]
Hence $\ker f$ is an extension of $\ker (f\circ i)$ and $\ker \partial
$, both of which are in $\cat{D}$ by Lemma~\ref{lem-extensions}.  So
$\ker f\in \cat{D}_{1}$.  Similarly, $\cok f=\cok \partial $, which is
in $\cat{D}_{1}$ by Lemma~\ref{lem-extensions}.   
\end{proof}

\begin{corollary}\label{cor-zaftig-generated-by}
Suppose $\cat{A}$ is an abelian category, and $\cat{D}$ is a full
subcategory of $\cat{A}$ closed under kernels and cokernels.  Let
$\cat{D}_{0}=\cat{D}$, and, for $n\geq 1$, define $\cat{D}_{n}$ to be
the full subcategory of extensions of $\cat{D}_{n-1}$.  Let
$\cat{E}=\bigcup _{n=0}^{\infty } \cat{D}_{n}$.  Then $\cat{E}$ is the
wide subcategory generated by $\cat{D}$.  
\end{corollary}

\begin{proof}
Note that the union that defines $\cat{E}$ is an increasing one, in the
sense that $\cat{D}_{n}\subseteq \cat{D}_{n+1}$.  This makes it clear
that $\cat{E}$ is closed under extensions.
Proposition~\ref{prop-extensions} implies that $\cat{E}$ is closed under
kernels and cokernels.  Therefore $\cat{E}$ is a wide subcategory.
Since any wide subcategory containing $\cat{D}$ must contain
each $\cat{D}_{n}$, the corollary follows.  
\end{proof}

\begin{theorem}\label{thm-modding-out}
Suppose $R$ is a coherent ring, and $\ideal{a}$ is a two-sided ideal of
$R$ that is finitely generated as a left ideal.  Then the map
$u\mathcolon \wide (R/\ideal{a})\xrightarrow{}\wide (R)$ is an
embedding.
\end{theorem}

\begin{proof}
It suffices to show that $gf(\cat{D})=\cat{D}$, where $\cat{D}$ is a
wide subcategory of $\cat{C}_{0}(R/\ideal{a})$.  According to
Corollary~\ref{cor-zaftig-generated-by}, $u(\cat{D})=\bigcup_{n}
\cat{D}_{n}$, where $\cat{D}_{n}$ is the collection of extensions of
$\cat{D}_{n-1}$.  We prove by induction on $n$ that if $\ideal{a}$ acts
trivially on some $M\in \cat{D}_{n}$, then in fact $M\in \cat{D}$.  The
base case of the induction is clear, since $\cat{D}_{0}=\cat{D}$.  Now
suppose our claim is true for $n-1$, and $\ideal{a}$ acts trivially on
$M\in \cat{D}_{n}$.  Write $M$ as an extension
\[
0 \xrightarrow{} M' \xrightarrow{} M \xrightarrow{} M'' \xrightarrow{} 0
\]
where $M',M''\in \cat{D}_{n-1}$.  Then $\ideal{a}$ acts trivially on
$M'$ and $M''$, so $M',M''\in \cat{D}$.  Furthermore, this is an
extension of $R/\ideal{a}$-modules; since $\cat{D}$ is a wide
subcategory of $R/\ideal{a}$-modules, $M\in \cat{D}$.  The induction is
complete, and we find that $gf(\cat{D})=\cat{D}$.
\end{proof}

\begin{theorem}\label{thm-final}
Suppose $R$ is a commutative coherent ring such that 
\[
f_{R}\mathcolon \wide (R)\xrightarrow{} \thick (\cat{D}(R))
\]
is an isomorphism, and $\ideal{a}$ is a finitely generated ideal of $R$.
Then $f_{R/\ideal{a}}$ is also an isomorphism.  In particular, $f_{R}$
is an isomorphism for all rings $R$ that are quotients of regular
commutative coherent rings by finitely generated ideals.  
\end{theorem}

\begin{proof}
The second statement follows immediately from the first and
Theorem~\ref{thm-classification}.  To prove the first statement, note
that, by hypothesis and Thomason's theorem~\ref{thm-thomason}, the map
\[
G_{R}=\alpha ji \mathcolon \widetilde{J}(\Spec R)\xrightarrow{}\wide
(R)
\]
is an isomorphism.  It suffices to show that $G_{R/\ideal{a}}$ is a
surjection.  So suppose $\cat{D}$ is a wide subcategory of finitely
presented $R/\ideal{a}$-modules.  Then $u\cat{D}=G(S)$ for some $S\in
\widetilde{J}(\Spec R)$.  Since $M_{(\ideal{p})}=0$ for any $R$-module
$M$ such that $\ideal{a}M=0$ and any $\ideal{p}$ not containing
$\ideal{a}$, we have $S=T\cup D(\ideal{a})$ for a unique $T\subseteq
V(\ideal{a})=\Spec (R/\ideal{a})$.  One can easily see that $T\in
J(\Spec R/\ideal{a})$, but we claim that in fact $T\in
\widetilde{J}(\Spec R/\ideal{a})$.  Indeed, $T=S\cap V(\ideal{a})$, so
this claim boils down to showing that the inclusion
$V(\ideal{a})\subseteq \Spec R$ is a proper map.  This is well-known;
the inclusion of any closed subset in any topological space is proper.

Naturally, we claim that $G_{R/\ideal{a}}(T)=\cat{D}$.  We have
$\cat{D}\subseteq u\cat{D}=G(S)\subseteq G(T)$, since $S\supseteq T$.
To show the converse, it suffices to show that
$uG_{R/\ideal{a}}(T)\subseteq u\cat{D}=G(S)$, since $u$ is an embedding by
Theorem~\ref{thm-modding-out}.  But any module $M$ in
$G_{R/\ideal{a}}(T)$ has $M_{(\ideal{p})}=0$ for all $\ideal{p}\in T$
and for all $\ideal{p}\in D(\ideal{a})$.  Therefore
$uG_{R/\ideal{a}}(T)\subseteq G(T\cup D(\ideal{a}))=G(S)$.  
\end{proof}

\begin{corollary}\label{cor-k}
Suppose $R$ is a finitely generated commutative $k$-algebra, where $k$
is a principal ideal domain.  Then every wide subcategory of finitely
generated $R$-modules is a Serre class.  
\end{corollary}

Indeed, any such $R$ is a quotient of a finitely generated polynomial
ring over $k$, which is regular by Hilbert's syzygy theorem, by a
finitely generated ideal.  This corollary covers most of the Noetherian
rings in common use, though of course it does not cover all of them.  We
remain convinced that $f$ should be an isomorphism for all commutative
coherent rings.

\section{Localizing subcategories}\label{sec-loc}

In this section, we relate certain subcategories of $R\Mod $ to
localizing subcategories of $\cat{D}(R)$.  Recall that a localizing
subcategory is a thick subcategory closed under arbitrary direct sums.
We use the known classification of localizing subcategories of
$\cat{D}(R)$ when $R$ is a Noetherian commutative ring to deduce a
classification of wide subcategories of $R$-modules closed under
arbitrat coproducts.  We know of no direct proof of this classification.

Let $\localizing (R)$ denote the lattice of wide subcategories of
$R$-modules closed under arbitrary coproducts, and let $\mloc
(\cat{D}(R))$ denote the lattice of localizing subcategories of
$\cat{D}(R)$.  Just as before, we can define a map $f\mathcolon
\localizing (R)\xrightarrow{}\mloc (R)$, where $f(\cat{C})$ is the
collection of all $X$ such that $H_{n}X\in \cat{C}$ for all $n$.  The
proof of Proposition~\ref{prop-zaftig-gives-thick} goes through without
difficulty to show that $f(\cat{C})$ is localizing.

Similarly, we can define $g\mathcolon \mloc
(\cat{D}(R))\xrightarrow{}\localizing (R)$ by letting $g(\cat{D})$ be
the smallest wide subcategory closed under coproducts containing all the
$H_{n}X$, for $X\in \cat{D}$ and for all $n$.  The proof of
Proposition~\ref{prop-adjoint} goes through without change, showing that
$g$ is left adjoint to $f$.

\begin{lemma}\label{lem-local}
For any ring $R$ and any wide subcategory $\cat{C}$ of $R$-modules
closed under coproducts, we have $gf(\cat{C})=\cat{C}$.
\end{lemma}

\begin{proof}
Since $g$ is left adjoint to $f$, $gf(\cat{C})\subseteq \cat{C}$.  But,
given $M\in \cat{C}$, $S^{0}M$ is in $f(\cat{C})$, and hence
$M=H_{0}S^{0}M \in gf(\cat{C})$.
\end{proof}

It would be surprising if an arbitrary localizing subcategory of
$\cat{D}(R)$ were determined by the homology groups of objects in
it, but this is nevertheless the case when $R$ is Noetherian and
commutative.  Given a prime ideal $\ideal{p}$ of such an $R$, denote by
$k_{\ideal{p}}$ the residue field $R_{(\ideal{p})}/\ideal{p}$ of
$\ideal{p}$.  

\begin{theorem}\label{thm-local}
Suppose $R$ is a Noetherian commutative ring.  Then $f\mathcolon
\localizing (R)\xrightarrow{}\mloc (\cat{D}(R))$ is an isomorphism.
Furthermore, there is an isomorphism between the Boolean algebra
$2^{\Spec R}$ and $\localizing (R)$ that takes a set $A$ of prime ideals
to the wide subcategory closed under coproducts generated by the
$k_{\ideal{p}}$ for $\ideal{p}\in A$.
\end{theorem}

\begin{proof}
We have a map $\alpha \mathcolon 2^{\Spec R}\xrightarrow{}\localizing
(R)$, defined as in the statement of the theorem.  The composition
$f\alpha $ is proved to be an isomorphism, for $R$ Noetherian and
commutative, in~\cite{neeman-derived} (see also~\cite[Sections 6 and
9]{hovey-axiomatic}).  Since $f$ is injective, we conclude that $f$, and
hence also $\alpha $, is an isomorphism.
\end{proof}

For example, the wide subcategory closed under coproducts of abelian
groups corresponding to the prime ideal $0$ is the collection of
rational vector spaces; the wide subcategory closer under coproducts
corresponding to the set $\{0,p \}$ is the collection of $p$-local
abelian groups.


\providecommand{\bysame}{\leavevmode\hbox to3em{\hrulefill}\thinspace}

\end{document}